\documentclass[12pt]{article}
\usepackage[a4paper,margin=1.0in]{geometry}
\usepackage[colorlinks,citecolor=magenta,linkcolor=black]{hyperref}
\pdfpagewidth=\paperwidth \pdfpageheight=\paperheight
\usepackage{amsfonts,amssymb,amsthm,amsmath,eucal,tabu,url}
\usepackage{pgf}
 \usepackage{array}
 \usepackage{tikz-cd}
 \usepackage{pstricks}
 \usepackage{pstricks-add}
 \usepackage{pgf,tikz}
 \usetikzlibrary{automata}
 \usetikzlibrary{arrows}
 \usepackage{indentfirst}
\usepackage{tabularx} 

\theoremstyle{plain}
\newtheorem{thm}{Theorem}

\newtheorem*{theoremA}{Main Theorem}

\newtheorem{conjecture}[thm]{Conjecture}

\theoremstyle{definition}
\newtheorem{definition}[thm]{Definition}
\newtheorem{remark}[thm]{Remark}

\newtheorem{thevarthm}[thm]{\varthmname}

\newenvironment{varthm*}[1]{\trivlist\item[]{\bf #1.}\it}{\endtrivlist}

\renewcommand\geq{\geqslant}

\renewcommand\leq{\leqslant}

\newcommand\be{\begin{eqnarray*}}
\newcommand\ee{\end{eqnarray*}}

\newcommand\newop[2]{\def#1{\mathop{\rm #2}\nolimits}}
\newop\log{log}
\newop\ord{ord}
\newop\Gal{Gal}
\newop\SL{SL}
\newop\Bl{Bl}
\newop\mult{mult}
\newop\mass{mass}
\newop\div{div}
\newop\codim{codim}
\newop\sing{sing}
\newop\vdim{vdim}
\newop\edim{edim}
\newop\Ass{Ass}
\newop\size{size}
\newop\reg{reg}
\newop\satdeg{satdeg}
\newop\supp{supp}
\newop\Neg{Neg}
\newop\Nef{Nef}
\newop\Nefh{Nef_H}
\newop\Eff{Eff}
\newop\Zar{Zar}
\newop\MB{MB}
\newop\MBxC{MB\mathit{(x,C)}}
\newop\NnB{NnB}
\newop\Bigg{Big}
\newop\Effbar{\overline{\Eff}}

\def\keywordname{{\bfseries Keywords}}%
\def\keywords#1{\par\addvspace\medskipamount{\rightskip=0pt plus1cm
\def\and{\ifhmode\unskip\nobreak\fi\ $\cdot$
}\noindent\keywordname\enspace\ignorespaces#1\par}}
\def\subclassname{{\bfseries Mathematics Subject Classification
(2020)}\enspace}
\def\subclass#1{\par\addvspace\medskipamount{\rightskip=0pt plus1cm
\def\and{\ifhmode\unskip\nobreak\fi\ $\cdot$
}\noindent\subclassname\ignorespaces#1\par}}

\begin{document}
\title{Weak Ziegler pairs of conic-line arrangements with ordinary singularities}
\author{Magdalena Lampa-Baczy\'nska and Daniel W\'ojcik}
\date{\today}
\maketitle

\thispagestyle{empty}
\begin{abstract}
We find the first weak Ziegler pair in the class of conic-line arrangements with ordinary quasi-homogeneous singularities defined over the rational numbers.
\keywords{weak Terao pair, conic-line arrangements, plane curve singularities}
\subclass{14N25, 14H50, 32S25}
\end{abstract}

In the present note we focus on conic-line arrangements in the plane with quasi-homogeneous ordinary singularities from the perspective of  \textbf{weak Ziegler pairs}. The foundations of this article come from an active area of research devoted to the freeness and nearly freeness of curve arrangements in the complex projective plane and the so-called Numerical Terao's Conjecture. This conjecture boils down to a very fundamental problem in combinatorial algebraic geometry, namely whether the \textbf{weak combinatorics} of a given arrangement determines the freeness. Let us recall the following definition.
\begin{definition}
Let $C = \{C_{1}, ..., C_{k}\} \subset \mathbb{P}^{2}_{\mathbb{C}}$ be a reduced curve such that each irreducible component $C_{i}$ is \textbf{smooth}. The weak combinatorics of $C$ is a vector of the form $(d_{1}, ..., d_{s}; q_{1}, ..., q_{p})$, where $d_{i}$ denotes the number of irreducible components of $C$ of degree $i$, and $q_{j}$ denotes the number of singular points of a curve $C$ of a given analytic type $Q_{j}$.
\end{definition}

Since we are working with conic-line arrangements $\mathcal{CL} \subset \mathbb{P}^{2}_{\mathbb{C}}$ with ordinary singularities, then the weak-combinatorics of such $\mathcal{CL}$ is the vector of the form $(d, k; t_{2}, ..., t_{d})$, where $t_{j}$ denotes the number of $j$-fold intersection points of $\mathcal{CL}$, i.e., points where exactly $j$ curves meet, $d\geq 1$ is the number of lines, and $k\geq 1$ is the number of smooth conics.

Having the above definition in hand, we can formulate the main open problem.
\begin{conjecture}[Numerical Terao's Conjecture]
Let $C_{1}, C_{2}$ be two reduced curves in $\mathbb{P}^{2}_{\mathbb{C}}$ such that their all irreducible components are smooth. Suppose that $C_{1}$ and $C_{2}$ have the same weak combinatorics and all singularities that our curves admit are quasi-homogeneous. Assume that $C_{1}$ is free, then $C_{2}$ has to be free.
\end{conjecture}
This conjecture, which is adjust to the world of curve arrangements with quasi-homogeneous singularities, is a natural generalization of the (classical) Terao's freeness conjecture formulated for central hyperplane arrangements. Recall that Terao's freeness conjecture predicts that the freeness property is determined by the intersection lattices. In order to understand the Numerical Terao's Conjecture, we need to construct examples of reduced curve arrangements with quasi-homogeneous singularities such that they have the same weak combinatorics, but different homological properties of the associated Milnor algebras. Before we explain in detail our result, we need some preparation.

Let $S := \mathbb{C}[x,y,z]$ denote the coordinate ring of $\mathbb{P}^{2}_{\mathbb{C}}$, and for a homogeneous polynomial $f \in S$ let $J_{f}$ denote the Jacobian ideal associated with $f$, i.e., the ideal of the form $J_{f} = \langle \partial_{x}\, f, \partial_{y} \, f, \partial_{z} \, f \rangle$.
\begin{definition}
\label{hom}
Let $C : f=0$ be a reduced curve in $\mathbb{P}^{2}_{\mathbb{C}}$ of degree $d$ given by $f \in S$. Denote by $M(f) := S/ J_{f}$ the Milnor algebra. We say that $C$ is $m$-syzygy when $M(f)$ has the following minimal graded free resolution:
$$0 \rightarrow \bigoplus_{i=1}^{m-2}S(-e_{i}) \rightarrow \bigoplus_{i=1}^{m}S(1-d - d_{i}) \rightarrow S^{3}(1-d)\rightarrow S$$
with $e_{1} \leq e_{2} \leq ... \leq e_{m-2}$ and $1\leq d_{1} \leq ... \leq d_{m}$. 
\end{definition}
\begin{definition}
The $m$-tuple $(d_{1}, ..., d_{m})$ in Definition \ref{hom} is called the exponents of $C$.
\end{definition}
In the present note we will work with a specific class of curves, namely nearly free curves \cite{DimcaSticlaru}.
\begin{definition}
A reduced plane curve $C \subset \mathbb{P}^{2}_{\mathbb{C}}$ of degree $d$ is nearly-free if $C$ is $3$-syzygy with $d_{1}+d_{2}=d$ and $d_{2}=d_{3}$.
\end{definition}
In order to define weak Ziegler pairs, we need the following definition.
\begin{definition}
Consider the graded $S$-module of Jacobian syzygies of $f$, namely $$AR(f)=\{(a,b,c)\in S^3 : af_x+bf_y+cf_z=0\}.$$
The minimal degree of non-trivial Jacobian relations for $f$ is defined to be 
$${\rm mdr}(f):=\min\{r : AR(f)_r\neq (0)\}.$$ 
In the light of Definition \ref{hom}, one has
$$d_{1} = {\rm mdr}(f).$$
\end{definition}
Finally, we can define the main object of our interests, following the lines of \cite[Definition 4.4]{Pokora1}.
\begin{definition}[Weak Ziegler pair]
Consider two reduced plane curves $C_{1}, C_{2} \subset \mathbb{P}^{2}_{\mathbb{C}}$ such that all irreducible components of $C_{1}$ and $C_{2}$ are smooth. We say that a pair $(C_{1},C_{2})$ forms a weak Ziegler pair if $C_{1}$ and $C_{2}$ have the same weak-combinatorics, but they have different minimal degrees of non-trivial Jacobian relations, i.e., ${\rm mdr}(C_{1}) \neq {\rm mdr}(C_{2})$.
\end{definition}

Having preparation done, we can formulate our main result devoted to the existence of the first weak Ziegler pair in the class of conic-line arrangements with ordinary quasi-homogeneous singularities.
\begin{theoremA}
There exists a pair of conic-line arrangements defined over the rational having the same weak-combinatorics of the form
$$(d,k;t_{2},t_{3},t_{4}) = (6,1; 6, 3, 2),$$
but having different minimal degrees of Jacobian relations. Furthermore, these two conic-line arrangements are nearly free.
\end{theoremA}

Our inspiration to look for such an example comes from a paper by Pokora \cite{Pokora}, where he studied free arrangements of lines and a conic with quasi-homogeneous ordinary singularities. We should also note here that our example gives a positive answer to \cite[Problem 4.8]{Pokora1}.

\begin{remark}
Let us recall that an ordinary singularity of multiplicity $m$ is quasi-homogeneous if $m<5$, see \cite[Exercise 7.31]{RCS}. It explains that our conic-line arrangements do indeed have quasi-homogeneous singularities. 
\end{remark}

Here is the proof of \textbf{Main Theorem}.
\begin{proof}
 We introduce conic-line arrangements $\mathcal{CL}_{1} \, : Q_{1} = 0$ and $\mathcal{CL}_{2} \, : Q_{2}=0$ by their defining polynomials, namely
 \begin{multline*}
 Q_{1}(x,y,z) = (x^2 +y^2 -25z^{2})(x-4z)(x+4z)(3x+4y)(x+3z)(6x+y+21z)(y+3z)
 \end{multline*}
and
\begin{multline*}
Q_{2}(x,y,z) = (x^2 +y^2 -25z^{2})(x+4z)x(x-4z)(y-3z)(y+3z)(3x-4y).
\end{multline*}
It is easy to check that both arrangements have the same weak combinatorics, namely
$$(d,k;t_{2},t_{3},t_{4},t_{5},t_{6}) = (6,1; 6, 3, 2).$$
Using \verb}SINGULAR} \cite{Singular}, we can compute the minimal free resolutions of the associated Milnor algebras. First of all, we have
$$(\mathcal{CL}_{1}) \, : \quad 0 \rightarrow S(-12) \rightarrow S(-11)^{3} \rightarrow S(-7)^{3} \rightarrow S,$$
which means that ${\rm mdr}(Q_{1}) = 4$, and our arrangement is nearly free with the exponents $(d_{1},d_{2},d_{3})=(4,4,4)$. 

In the case of the second arrangement, we have
$$(\mathcal{CL}_{2}) \, : \quad 0 \rightarrow S(-13) \rightarrow S(-12)^2 \oplus S(-10) \rightarrow S(-7)^{3} \rightarrow S,$$
which means that ${\rm mdr}(Q_{2}) = 3$, and our arrangement is nearly free with the exponents $(d_{1},d_{2},d_{3})=(3,5,5)$. This completes the proof.
\end{proof}
To the best of our knowledge, our example is the first of its kind, since our weak Ziegler pair is in the class of nearly free arrangements with quasi-homogeneous singularities. So far we have seen examples of pairs where the homological descriptions of the curves are different, e.g. one curve is free, the second one is nearly free. 

For the completeness of our presentation, and in order to give some feeling about the construction, we deliver below figures of our conic-line arrangements.

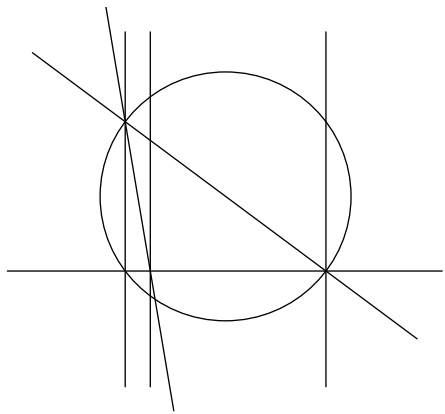
\begin{figure}[h]
\centering
\begin{tikzpicture}[line cap=round,line join=round,>=triangle 45,x=1cm,y=1cm,scale=0.33]
\clip(-12.697691760870853,-10.652515046699648) rectangle (15.63332662775391,7.609639663737902);
\draw [line width=0.5pt] (0,0) circle (5cm);
\draw [line width=0.5pt] (4,-7.652515046699648) -- (4,6.609639663737902);
\draw [line width=0.5pt] (-4,-7.652515046699648) -- (-4,6.609639663737902);
\draw [line width=0.5pt,domain=-7.697691760870853:7.63332662775391] plot(\x,{(-0-3*\x)/4});
\draw [line width=0.5pt] (-3,-7.652515046699648) -- (-3,6.609639663737902);
\draw [line width=0.5pt,domain=-12.697691760870853:-2.063332662775391] plot(\x,{(-21-6*\x)/1});
\draw [line width=0.5pt,domain=-8.697691760870853:8.63332662775391] plot(\x,{(-3-0*\x)/1});
\end{tikzpicture}
\caption{Geometric presentation of arrangement $\mathcal{CL}_{1}$.}
\end{figure}
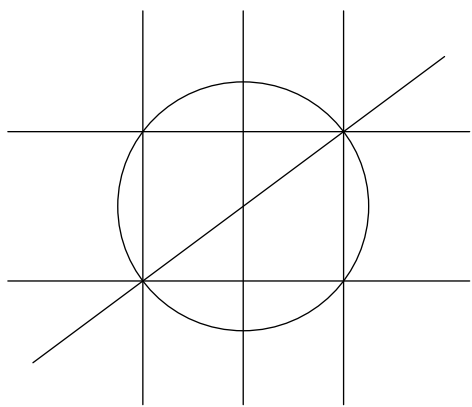
\begin{figure}[h]
\centering
    \begin{tikzpicture}[line cap=round,line join=round,>=triangle 45,x=1cm,y=1cm,scale=0.33]
\clip(-17.380556201367607,-10.969801325407138) rectangle (18.024648139313477,11.85236871997969);
\draw [line width=0.5pt] (0,0) circle (5cm);
\draw [line width=0.5pt] (4,-7.969801325407138) -- (4,7.85236871997969);
\draw [line width=0.5pt] (-4,-7.969801325407138) -- (-4,7.85236871997969);
\draw [line width=0.5pt] (0,-7.969801325407138) -- (0,7.85236871997969);
\draw [line width=0.5pt,domain=-9.380556201367607:9.024648139313477] plot(\x,{(-3-0*\x)/1});
\draw [line width=0.5pt,domain=-9.380556201367607:9.024648139313477] plot(\x,{(--3-0*\x)/1});
\draw [line width=0.5pt,domain=-8.380556201367607:8.024648139313477] plot(\x,{(-0-3*\x)/-4});
\end{tikzpicture}
\caption{Geometric presentation of arrangement $\mathcal{CL}_{2}$.}
\end{figure}
\section*{Data availability statement}
I do not analyse or generate any datasets, because this work proceeds within a theoretical and mathematical approach. 
\section*{Acknowledgement}
We want to thank Piotr Pokora for his suggestions and remarks about the content of the note.

\vskip 0.5 cm
\bigskip
Magdalena Lampa-Baczy\'nska,
Department of Mathematics,
University of the National Education Commission Krakow,
Podchor\c a\.zych 2,
PL-30-084 Krak\'ow, Poland. \\
\nopagebreak
\textit{E-mail address:} \texttt{magdalena.lampa-baczynska@up.krakow.pl}
\bigskip

Daniel W\'ojcik,
Department of Mathematics,
University of the National Education Commission Krakow,
Podchor\c a\.zych 2,
PL-30-084 Krak\'ow, Poland. \\
\nopagebreak
\textit{E-mail address:} \texttt{daniel.wojcik@up.krakow.pl}

\begin{thebibliography}{000}
\bibitem{Singular}
W.~Decker, G.-M. Greuel, G.~Pfister, and H.~Sch\"onemann,
\newblock {\sc Singular} {4-1-1} --- {A} computer algebra system for polynomial computations. \newblock \url{http://www.singular.uni-kl.de}, 2018.

\bibitem{RCS}  
A. Dimca,  \textit{Topics on real and complex singularities. An introduction.} Advanced Lectures in Mathematics. Braunschweig/Wiesbaden: Friedr. Vieweg \& Sohn. 1987.

\bibitem{DimcaSticlaru}
A. Dimca, G. Sticlaru, Free and Nearly Free Curves vs. Rational Cuspidal Plane Curves. \textit{Publ. Res. Inst. Math. Sci.} \textbf{54(1)}: 163 -- 179 (2018). 

\bibitem{Pokora}
P. Pokora, Freeness of arrangements of lines and a conic with ordinary quasi-homogeneous singularities. \textbf{arXiv:2312.13052}.

\bibitem{Pokora1}
P. Pokora, Singular plane curves and their interplay with
algebra and combinatorics. \textbf{arXiv:2403.13377}.
\end{thebibliography}
\end{document}